\documentclass{article}

\usepackage{referencing}
\usepackage{authblk}
\usepackage{amsmath, amsthm, amssymb}
\usepackage{graphicx}
\usepackage{xcolor}
\usepackage{tikz}
\usetikzlibrary{shapes.geometric}
\tikzstyle{V} = [circle, draw, thin, fill=white, scale=0.5]
\tikzstyle{E} = [->, very thick]
\usepackage{setspace}
\usepackage{thmtools,thm-restate}

\usepackage{multicol}
\usepackage{mathrsfs}
\usepackage{abstract}
\usepackage{enumitem}
\title{Residual Finiteness and Related Properties in Monounary Algebras and their Direct Products}
\author{Bill de Witt}
\affil{Mathematical Institute, School of Mathematics and Statistics, University of St Andrews, St Andrews, Fife, KY16 9SS, UK. \textit{Email}: bldw@st-andrews.ac.uk}
\date{}

\begin{document}

\maketitle

\begin{abstract}
In this paper we discuss the relationship between direct products of monounary algebras and their components, with respect to the properties of residual finiteness, strong/weak subalgebra separability, and complete separability. For each of these properties \(\mathcal{P}\), we give a criterion \(\mathcal{C_P}\) such that a monounary algebra \(A\) has property \(\mathcal{P}\) if and only if it satisfies \(\mathcal{C_P}\). We also show that for a direct product \(A\times B\) of monounary algebras, \(A\times B\) has property \(\mathcal{P}\) if and only if one of the following is true: either both \(A\) and \(B\) have property \(\mathcal{P}\), or at least one of \(A\) or \(B\) are \textit{backwards-bounded}, a special property which dominates direct products and which guarantees all \(\mathcal{P}\) hold.

\end{abstract}

\section{Summary of Results}

Monounary algebras are the simplest types of algebraic structure which are not entirely trivial, and yet display some interesting structure and behaviours. In this paper we work with residual finiteness and the related properties of strong/weak subalgebra separability and complete separability. Throughout this paper we will use \(\mathbb{N}\) to denote the set of non-negative integers, and \(\mathbb{N}^+\) to denote the set of positive integers. For a monounary algebra \((A,f)\) we use, for \(x\in A\) and \(n\in\mathbb{N}\) the notion of preimage sets \(f^{-n}(x)\) (defined in Section~\ref{Basics}) to give necessary and sufficient conditions for \(A\) to have these properties (Theorems \ref{MainRFThm}, \ref{SSEquivalence}, and \ref{CStheorem}). Specifically:

\begin{description}
    \item [Residual Finiteness]: For all \(x,y\in A\) such that \(x\neq y\) and \(f(x)=f(y)\), there exists \(n\in\mathbb{N}\) such that either \(f^{-n}(x)=\varnothing\) or \(f^{-n}(y)=\varnothing\).
    \item [Strong/Weak Subalgebra Separability]: For all \(x\in A\), either there exists \(n\in\mathbb{N}\) such that \(f^{-n}(x)=\varnothing\), or \(x\) is in a cycle.
    \item [Complete Separability]: For all \(a\in A\) there exists \(n\in\mathbb{N}\) such that \(f^{-n}(a)\setminus\bigcup_{i=0}^{n-1}f^{-i}(a)=\varnothing\).
\end{description}

We then consider direct products, and take into particular consideration algebras where for every \(x\in A\) there exists \(n\in\mathbb{N}\) such that \(f^{-n}(x)=\varnothing\), which we call \textit{backwards-bounded}. We show in Theorems \ref{RFproducttheorem}, \ref{SSproducttheorem}, and \ref{CSproducttheorem} that when it comes to direct products, all of these properties behave in the same way. More precisely, we show that a direct product has property \(\mathcal{P}\) if and only if one of the following is true: both components have property \(\mathcal{P}\), or at least one of them is backwards-bounded.

\section{Basics of Monounary Algebras and Residual Finiteness} \label{Basics}
We will begin with an introduction to monounary algebras, and prove some results regarding their structure which are both useful in the overall context of this paper, but also are helpful for the reader to visualise these algebraic structures. A general overview of monounary algebras can be found in \cite{jakubikova2009monounary}

A \textit{unary operation} is a function from a set to itself, and a \textit{monounary algebra} is a set together with a single unary operation defined on it. Note that we will usually identify a monounary algebra with its underlying set, and as such will omit mentioning the function where it is not necessary. A monounary algebra \((A,f)\) can be visualised in a natural way, as a directed graph; the vertices are the elements of \(A\), and for all \(a\in A\) there is a directed edge from \(a\) to its image \(f(a)\). Note that there is exactly one out-edge at each vertex. We now define a few specific monounary algebras, which will be used in the paper:

\begin{figure}
    \centering
    \begin{tikzpicture}[scale=2, >=stealth]
        \foreach \x/\place in {{1/(0,0)},{2/(0.5,0.5)},{3/(0,1)},{4/(-0.5,0.5)}}
            \node[V] (a\x) at \place {};
            \draw[E] (a1) -- (a2);
            \draw[E] (a2) -- (a3);
            \draw[E] (a3) -- (a4);
            \draw[E] (a4) -- (a1);
        \foreach \x/\place in {{1/(-1.5,0.5)},{2/(-2,0.5)},{3/(-2.5,0.5)},{4/(-3,0.5)}}
            \node[V] (b\x) at \place {};
            \draw[E] (b1) to[out=-45,in=180] (-1.2,0.25) to[out=0,in=0] (-1.2,0.75) to[out=180,in=45] (b1);
            \draw[E] (b2) -- (b1);
            \draw[E] (b3) -- (b2);
            \draw[E] (b4) -- (b3);
        \foreach \n/\x/\y in {{1/1.5/0.3},{2/2/0.3},{3/2.5/0.3},{4/3/0.3}}
            {\node[V] (d\n) at (\x,\y) {};
            \draw[E] (d\n) to[out=135,in=-90] (\x-0.2,\y+0.25) to[out=90,in=90] (\x+0.2,\y+0.25) to[out=-90,in=45] (d\n);};
        \draw (-0.75,1.25) rectangle (0.75,-0.25);
        \draw (-3.25,1.25) rectangle (-1,-0.25);
        \draw (1,1.25) rectangle (3.5,-0.25);
        \node (L4) at (-3,1) {\(L_4\)};
        \node (C4) at (-0.5,1) {\(C_4\)};
        \node (T4) at (1.25,1) {\(T_4\)};
         \foreach \n/\place in {{-1/(-0.5,-1)},{0/(0,-1)},{1/(0.5,-1)}}
            \node[V] (z\n) at \place {};
        \draw[E] (z-1) -- (z0);
        \draw[E] (z0) -- (z1);
        \draw[dotted,thick] (-1,-1) -- (z-1);
        \draw[dotted,thick] (z1) -- (1,-1);
        \draw (-1.25,-0.5) rectangle (1.25,-1.5);
        \node (Z) at (-1,-0.7) {\(\mathbb{Z}\)};
    \end{tikzpicture}
    \caption{The monounary algebras $L_4$, $C_4$, $T_4$ and $\mathbb{Z}$.}
    \label{fig:C4}
\end{figure}
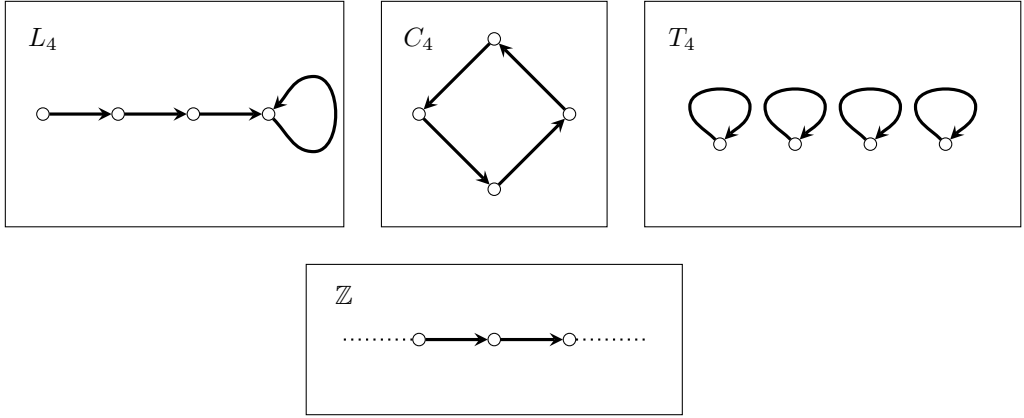

\begin{example} \label{exampleMUAs}
\begin{enumerate}
    \item The \(n\)-\textit{line} is the monounary algebra \(L_n=(\{0,1,\dots,n-1\},x\mapsto\max\{0,x-1\})\).
    \item The \(n\)-\textit{cycle} is the monounary algebra \(C_n=(\{0,1,\dots,n-1\},x\mapsto (x+1)\) mod \(n)\).
    \item The \(n\)-\textit{trivial} monounary algebra on \(n\) points is \(T_n=(\{0,1,\dots,n-1\},x\mapsto x)\).
    \item The \textit{bi-infinite path} \((\mathbb{Z},x\mapsto x+1)\), hereafter referred to as \(\mathbb{Z}\).
\end{enumerate}
\end{example}

Note that \(L_n, C_n\) and \(T_n\) are defined for all \(n\in\mathbb{N}^+\).

Some specific instances of these algebras are depicted in Figure \ref{fig:C4}. 

We will be using a number of results about the structure of these graphs, rephrased as results about monounary algebras. Many of these results can be found in \cite{ganyushkin2008classical} in more detail. Relevant graph theory definitions and results can be found in \cite{west2001introduction}, and a discussion of the combinatorial aspects of monounary algebras can be found in \cite{katai2019fine}.

One particular graph theoretic property we will deal with is connectedness. A monounary algebra is called \textit{connected} if the corresponding undirected graph is connected.  For the purposes of this paper, it will almost always be sufficient to prove results for connected monounary algebras.

As is the usual for algebraic structures, we have notions of subalgebras and homomorphisms. For a monounary algebra \((A,f)\) a \textit{subalgebra} \((S,g)\) of \(A\) is a subset \(S\subseteq A\) such that \(f(S)\subseteq S\) with a unary operation \(g=f|_S\), the restriction of \(f\) to \(S\). A function \(\phi:A\to B\) between two monounary algebras \((A,f_1)\) and \((B,f_2)\) is a \textit{homomorphism} if \(\phi(f_1(a))=f_2(\phi(a))\) for all \(a\in A\). An \textit{isomorphism} is a homomorphism which is a bijection.

\begin{lemma} \label{UniqueCyclesLemma}
Let \(A\) be a  monounary algebra. Then:
\begin{enumerate}
    \item If \(A\) is finite and non-empty, there is a subalgebra \(C\leq A\) which is a cycle (i.e. \(C\cong C_k\) for some \(k\in\mathbb{N}^+\)).
    \item If \(A\) is connected and there exists a subalgebra \(C\) isomorphic to a cycle, then it is the unique such subalgebra, and is contained in every non-empty subalgebra of \(A\).
\end{enumerate} 
\end{lemma}

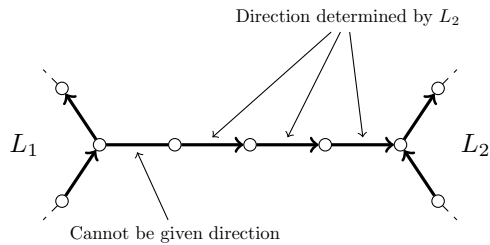
\begin{figure}
    \centering
    \begin{tikzpicture}
        \foreach \n/\place in {{1/(-2,0)},{2/(-1,0)},{3/(0,0)},{4/(1,0)},{5/(2,0)},{6/(-2.5,0.75)},{7/(-2.5,-0.75)},{8/(2.5,0.75)},{9/(2.5,-0.75)}}
            \node[V] (p\n) at \place {};
        \draw[very thick] (p1) -- (p2);
        \draw[E] (p2) -- (p3);
        \draw[E] (p3) -- (p4);
        \draw[E] (p4) -- (p5);
        \draw[E] (p5) -- (p8);
        \draw[E] (p7) -- (p1);
        \draw[E] (p9) -- (p5);
        \draw[E] (p1) -- (p6);
        \draw[dashed] (p6) -- (-2.75,1);
        \draw[dashed] (p7) -- (-2.75,-1);
        \draw[dashed] (p8) -- (2.75,1);
        \draw[dashed] (p9) -- (2.75,-1);
        \node (L1) at (-3,0) {\(L_1\)};
        \node (L2) at (3,0) {\(L_2\)};
        \node[scale=0.7, above] (detstat) at (1.3,1.5) {Direction determined by \(L_2\)};
        \draw[->] (detstat) -- (-0.5,0.1);
        \draw[->] (detstat) -- (0.5,0.1);
        \draw[->] (detstat) -- (1.5,0.1);
        \node[scale=0.7,below] (imp) at (-1,-1) {Cannot be given direction};
        \draw[->] (imp) to (-1.5,-0.1);
    \end{tikzpicture}
    \caption{There cannot be a path between two cycles.}
    \label{fig:looppath}
\end{figure}

We omit the proof of this Lemma, as it is sufficiently simple. The important observation to show uniqueness is that paths cannot come out of a cycle, they must go into cycles (see Figure \ref{fig:looppath}).

For a monounary algebra \((A,f)\), and a subset \(S\subseteq A\), we may want to consider the \textit{subalgebra generated by }\(S\), denoted \(\langle S\rangle\). By this we mean the smallest subalgebra of \(A\) containing \(S\), which can be constructed as follows:
\[\langle S\rangle = \{f^n(s):n\in\mathbb{N}\text{, }s\in S\}.\]

Residual finiteness is a property which has been studied in depth for a number of algebraic structures (groups in particular) for many decades. One can find discussions of residual finiteness in groups in \cite[Chapter 9]{MR0332990} as well as in \cite[Chapter 2]{ceccherini2010cellular}. The notion has also been investigated in general algebraic structures, such as in \cite{evans1969some}, which looks at the relationship between residual finiteness and other algebraic properties. A recent work in this area is \cite{mayr2018finiteness}, which discusses residual finiteness of direct products in congruence modular varieties. It is with the goal of understanding how the property is reflected in general algebraic structures, that we investigate its nature in this specific type of structure.

We now define the property of residual finiteness for monounary algebras.

\begin{definition}
A monounary algebra \(A\) is \textit{residually finite} if for all distinct \(x,y\in A\), there exists a finite monounary algebra \(F\) and a homomorphism \(\phi:A\to F\) such that \(\phi(x)\neq\phi(y)\).
\end{definition}

We give the specific example of \(\mathbb{Z}\), and show that it is residually finite. This will be used in the main classification theorem for residual finiteness in Section~\ref{sec:RFCrit}.

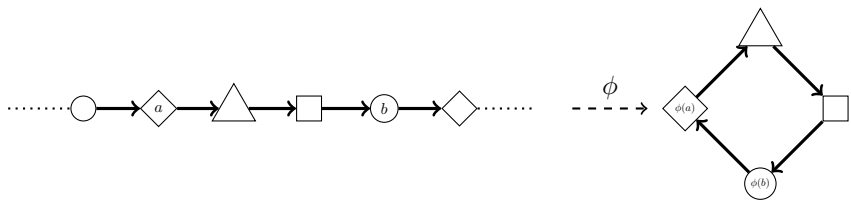
\begin{figure}
    \centering
    \begin{tikzpicture}
        \node[circle, draw, fill=none] (z0) at (-6,0) {};
        \node[diamond, draw, fill=none,scale=0.6] (z1) at (-5,0) {\(a\)};
        \node[regular polygon,regular polygon sides=3, draw] (z2) at (-4,0) {};
        \node[regular polygon,regular polygon sides=4, draw, fill=none] (z3) at (-3,0) {};
        \node[circle, draw, fill=none,scale=0.6] (z4) at (-2,0) {\(b\)};
        \node[diamond, draw, fill=none] (z5) at (-1,0) {};
        \draw[E] (z0) -- (z1);
        \draw[E] (z1) -- (z2);
        \draw[E] (z2) -- (z3);
        \draw[E] (z3) -- (z4);
        \draw[E] (z4) -- (z5);
        \draw[dotted, thick] (-7,0) -- (z0);
        \draw[dotted, thick] (z5) -- (0,0);
        \node[diamond, draw, fill=none,scale=0.4] (c1) at (2,0) {\(\phi(a)\)};
        \node[regular polygon,regular polygon sides=3, draw] (c2) at (3,1) {};
        \node[regular polygon,regular polygon sides=4, draw, fill=none] (c3) at (4,0) {};
        \node[circle, draw, fill=none,scale=0.4] (c4) at (3,-1) {\(\phi(b)\)};
        \draw[E] (c1) -- (c2);
        \draw[E] (c2) -- (c3);
        \draw[E] (c3) -- (c4);
        \draw[E] (c4) -- (c1);
        \draw[dashed,->,thick] (0.5,0) -- (1.5,0);
        \node[above] (phi) at (1,0) {\(\phi\)};
    \end{tikzpicture}
    \caption{The homomorphism $phi$ separates $a$ and $b$ in a finite algebra.}
    \label{fig:ZRF}
\end{figure}

\begin{lemma} \label{ZRF}
The monounary algebra \(\mathbb{Z}\) is residually finite.
\begin{proof}
Let \(a\),\(b\in\mathbb{Z}\) with \(a\neq b\),  and set \(m=|b-a|+1\). Then we construct a map \(\phi:\mathbb{Z}\to C_m\) defined by \(\phi(n) = n\pmod{m}\).
It is then easy to verify that this is a homomorphism, and that \(\phi(a)\neq\phi(b)\). An example is depicted in Figure \ref{fig:ZRF}.
\end{proof}
\end{lemma}

One of the main motivations for this paper is to study how residual finiteness interacts with direct products. The following is well-known in universal algebra:

\begin{lemma}
Let \(A\) and \(B\) be residually finite algebras of the same type. Then the direct product \(A\times B\) is residually finite.
\begin{proof}
If two pairs \((a_1,b_1),(a_2,b_2)\in A\times B\) are not equal, then they differ in at least one component. Assume without loss of generality that \(a_1\neq a_2\). Then as \(A\) is residually finite there exists a homomorphism \(\phi:A\to F\) where \(F\) is finite and \(\phi(a_1)\neq\phi(a_2)\). Thus, letting \(\pi_1:A\times B\to A\) be the projection onto the first co-ordinate, we have that \(\theta=\pi_1\circ\phi:A\times B\to F\) is a homomorphism such that \(\theta(a_1)\neq\theta(a_2)\).
\end{proof}
\end{lemma}

But the converse is more challenging. That is: if \(A\times B\) is residually finite, is it true that both \(A\) and \(B\) are residually finite? It was shown in \cite{gray2009residual} that this is true for many well-studied classes of algebras, via the following proposition.

\begin{proposition}
Let \(A\) and \(B\) be algebras, and suppose that \(A\) contains an idempotent. If \(A\times B\) is residually finite then \(B\) is residually finite.
\end{proposition}

Note: an idempotent is an element \(e\) such that \(f(e,e,\dots,e)=e\) for every operation \(f\). Equivalently, it is an element such that \(\{e\}\) is a subalgebra.

From this proposition it follows that the condition is true for any class which always contains idempotents, such as: groups, rings, monoids, semilattices, loops, quasirings etc.

Less obviously, in \cite{gray2009residual} the analogous assertion was also shown for the variety of semigroups, which do not necessarily have idempotents. However it is not true in unary algebras. Specifically, \cite{gray2009residual} provided an example of a residually finite product of monounary algebras, one  of which is not residually finite, as well as an example of two biunary algebras neither of which are residually finite, but whose product is. We will show in Section~\ref{sec:Prod}    that in monounary algebras, at least one component must be residually finite to obtain residual finiteness in the direct product.

\section{Preliminaries on Preimages}

We now introduce some notation for a key concept featured in the major results in this paper, that of preimage sets.

\begin{notation}
Let \((A,f)\) be a monounary algebra, and \(a\in A\). Then for \(n\in\mathbb{N}\) we define \[f^{-n}(a) = \{b\in A:f^n(b)=a\}.\]
In graphical terms, this would be the set of points from which there is a walk of length \(n\) terminating at \(a\).
\end{notation}

It turns out that this is the most important thing to consider when dealing residual finiteness of a monounary algebra, and so this chapter is dedicated to their properties, and some constructions we can use them for.

\begin{definition}
For a monounary algebra \(A\), we say a point \(a\in A\) is \textit{backwards eternal} if \(f^{-n}(a)\neq\varnothing\) for all \(n\in\mathbb{N}\).
\end{definition}

Backwards eternality is particularly useful, and so it is worth clarifying that it is not equivalent to a point being at the end of a path of infinite length, it can also be at the end of infinitely many finite paths with no upper bound on length, as in the following example.

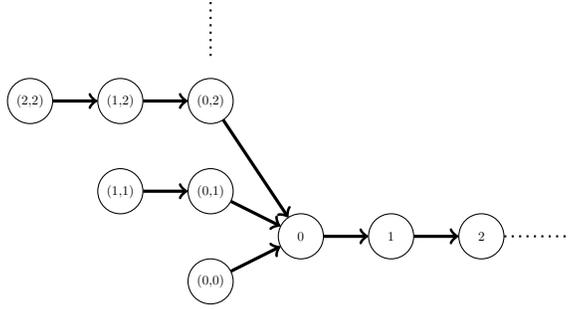
\begin{figure}
    \centering
    \begin{tikzpicture}[scale=1.2]
        \foreach \x/\place/\n in {{0/(0,0)/0},{1/(1,0)/1},{2/(2,0)/2},{00/(-1,-0.5)/(0,0)},{01/(-1,0.5)/(0,1)},{11/(-2,0.5)/(1,1)},{02/(-1,1.5)/(0,2)},{12/(-2,1.5)/(1,2)},{22/(-3,1.5)/(2,2)}}
            \node[V, minimum size=12mm] (a\x) at \place {\n};
        \foreach \o/\i in {{22/12},{12/02},{02/0},{11/01},{01/0},{00/0},{0/1},{1/2}}
            \draw[E] (a\o) -- (a\i);
        \draw[dotted, thick] (a2) -- (3,0);
        \draw[dotted, thick] (-1,2) -- (-1,2.6);
    \end{tikzpicture}
    \caption{A monounary algebra containing a backwards eternal point, but no backwards infinite paths.}
    \label{fig:Bi-eternalExample}
\end{figure}

\begin{example} \label{nonZbieternal}
Let \(A=\mathbb{N}\cup\{(a,b)\in\mathbb{N}^2:a\leq b\}\), and define a unary operation by 
\begin{equation*}
    f(x) =
    \begin{cases}
        x+1 & \text{if }x\in\mathbb{N},\\
        0 & \text{if }x\in\{0\}\times\mathbb{N},\\
        x-(1,0) &\text{otherwise}.\\
    \end{cases}
\end{equation*}
This monounary algebra is depicted in Figure \ref{fig:Bi-eternalExample}. Every point in \(\mathbb{N}\) is backwards eternal. 
\end{example}

Note that, for simplicity, in the following lemma we identify \(f^n(a)\) with the set \(\{f^n(a)\}\) for \(n\in\mathbb{N}\).
\begin{lemma} \label{PreimageLemmas}
 Let \((A,f)\) be a monounary algebra, and \(a,x\in A\). Then:
\begin{enumerate}
    \item [\textit{i)}] \(f^n(f^{-m}(a))\subseteq f^{n-m}(a)\) for all \(n\),\(m\in\mathbb{N}\);
    \item [\textit{ii)}] if there exists \(n\in\mathbb{N}\) such that \(f^{-n}(a)=\varnothing\), then \(f^{-m}(a)=\varnothing\) for all \(m\geq n\);
    \item [\textit{iii)}] if \(a\) is in a cycle then \(a\) is backwards eternal;
    \item [\textit{iv)}] if \(a\) is not in a cycle, then \(f^{-n}(a)\cap f^{-m}(a)=\varnothing\) for any distinct \(n\),\(m\in\mathbb{N}\);
    \item [\textit{v)}] \(f(x)\in f^{-n}(a)\) if and only if \(x\in f^{-(n+1)}(a)\).
\end{enumerate}
\begin{proof}
\textit{i)} First note that \(f^m(f^{-m}(a))=a\) by definition, and so the result follows trivially for \(n\geq m\).
For \(n<m\), note that for all \(b\in f^{-m}(a)\) we have that \(a=f^m(b)=f^{m-n}(f^n(b))\), so \(f^n(b)\in f^{n-m}(a)\).

\textit{ii)} Follows immediately from \textit{i)}.

\textit{iii)} Let \(k\) be the length of the cycle containing \(a\). Then \(a\in f^{-mk}(a)\) for all \(m\in\mathbb{N}\). So by \textit{ii)} we have that \(f^{-n}(a)\neq\varnothing\) for all \(n\in\mathbb{N}\). 

\textit{iv)} Assume not, and that \(n>m\). Then there exists \(b\in f^{-n}(a)\cap f^{-m}(a)\). So \(f^n(b)=f^m(b)=a\), and it follows that \(f^{n-m}(a)=a\), so \(a\) is in a cycle, a contradiction.

\textit{v)} Follows immediately from the definition.
\end{proof}
\end{lemma}

These results are sufficiently intuitive that they will be used without explicit reference in later sections. The following is slightly more complicated, and deals with the interactions between pre-images of distinct points.

\begin{lemma} \label{SeparableBackpaths}
Let \((A,f)\) be a connected monounary algebra, and let \(a,b\in A\) be two distinct elements. Then precisely one of the following is true:
\begin{enumerate}[label=(\arabic*)]
    \item  There exists \(n\in\mathbb{N}\) such that \(f^n(a)=b\) or \(f^n(b)=a\);
    \item  \((\bigcup_{n\in\mathbb{N}}f^{-n}(a))\cap(\bigcup_{n\in\mathbb{N}}f^{-n}(b))=\varnothing\).
\end{enumerate}

Additionally, in the case of \textit{(2)}, there exist \(n\),\(m\in\mathbb{N}\) such that \(f^n(a)=f^m(b)\).

\begin{proof}
It is clear that if \textit{(1)} is true then \textit{(2)} is false. It is then sufficient to show that if \textit{(1)} is false, then \textit{(2)} must be true.

We show the contrapositive. Assume \textit{(2)} is false. It therefore follows that there exists \(x\in A\) such that \(x\in f^{-n}(a)\cap f^{-m}(b)\) for some \(n\),\(m\in\mathbb{N}\). We thus have that \(f^n(x)=a\) and \(f^m(x)=b\). Assuming without loss of generality that \(n<m\), we get that \(b=f^m(x)=f^{m-n}(f^n(x))=f^{m-n}(a)\), so \textit{(1)} is true.

Finally, for the additional condition, when \textit{(2)} is true, as \(A\) is connected there exists an (undirected) path \((a=a_0,a_1,\dots,a_k=b)\). If this path were a directed path, then we would have \textit{(1)}, so this is not the case.

As an out-edge corresponds to the action of the function \(f\), there cannot be two out-edges at a vertex \(a_i\). Thus for a path to not be directed there must exist an \(a_i\) with two in-edges, or in other words, \(a_i=f(a_{i-1})=f(a_{i+1})\). Since \(f(a_{i-1})=a_i\neq a_{i-2}\) we have that \(f(a_{i-2})=a_{i-1}\), or in other words, the edge goes from \(a_{i-2}\) to \(a_{i-1}\). Repeating this process until we reach \(a\) gives us that \(f^i(a)=a_i\) and a similar process from \(a_{i+1}\) gives us that \(f^{k-i}(b)=a_i=f^i(a)\).
\end{proof}
\end{lemma}

\begin{figure}
    \centering
    \begin{tikzpicture} [scale=1.25]
        \foreach \x/\place/\n in {{1/(-1.5,0.75)/b},{2/(-2.5,0.75)/},{3/(-3.5,0.75)/},{4/(-4.5,0.75)/a}}
            \node[V,minimum size=20] (b\x) at \place {\n};
        \draw[E] (b2) -- (b1);
        \draw[->,dotted, thick] (b3) -- (b2);
        \draw[E] (b4) -- (b3);
        \foreach \x/\place/\n in {{1/(0,0)/a},{2/(2,0)/b},{3/(0.33,0.75)/},{4/(1.66,0.75)/},{5/(0.66,1.5)/},{6/(1.33,1.5)/},{7/(1,2.25)/}}
            \node[V,minimum size=20] (a\x) at \place {\n};
        \foreach \o/\i in {{1/3},{2/4},{5/7},{6/7}}
            \draw[E] (a\o) -- (a\i);
        \draw[dotted,->, thick] (a3)--(a5);
        \draw[dotted,->, thick] (a4)--(a6);
        \draw (a1) -- (-0.75,-1);
        \draw (a1) -- (0.75,-1);
        \draw (a2) -- (1.25,-1);
        \draw (a2) -- (2.75,-1);
        \node[scale=0.5] (ba) at (0,-0.8) {\(\displaystyle{\bigcup_{n\in\mathbb{N}}}f^{-n}(a)\)};
        \node[scale=0.5] (bb) at (2,-0.8) {\(\displaystyle{\bigcup_{n\in\mathbb{N}}}f^{-n}(b)\)};
    \end{tikzpicture}
    \caption{The two possibilities as discussed in Lemma \protect\ref{SeparableBackpaths}.}
    \label{fig:SeparableBackpaths}
\end{figure}
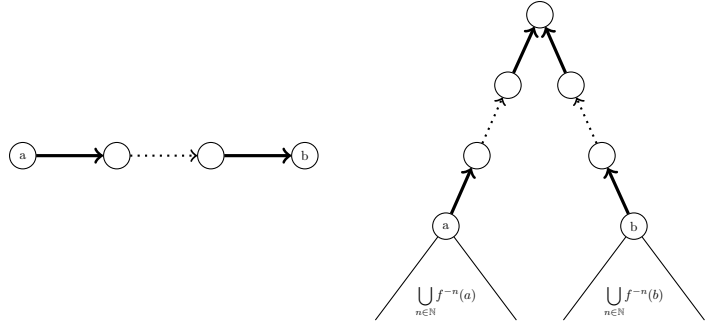

The two possibilities are depicted in Figure \ref{fig:SeparableBackpaths}.

We will now use preimage sets to construct some homomorphisms which will be used in later sections.

\begin{lemma} \label{cyclehomom}
Let \((A,f)\) be a connected monounary algebra containing a cycle of length \(k\). Then there exists a homomorphism \(\phi:A\to C_k\).
\begin{proof}
Fix a point \(a\in A\) contained in the cycle. Note that since \(A\) is connected, we have that for all \(y\in A\setminus\{a\}\) there exists \(n\in\mathbb{N}\) such that \(f^n(y)\) is in the cycle. Thus there exists a minimal \(n_y\in\mathbb{N}\) such that \(f^{n_y}(y)=a\). Then since the unary operation on \(C_k\) (which we will call \(f_{C_k}\)) is a bijection, \(f_{C_k}^{-n}\) is a well defined function for all \(n\in\mathbb{N}\). Thus we define a function \(\phi:A\to C_k\) as follows:
\begin{equation*}
    \phi(x) =
    \begin{cases}
        0 & \text{if }a=x,\\
        f_{C_k}^{-n_x}(0) &\text{otherwise}.
    \end{cases}
\end{equation*}
We show this is a homomorphism. If \(n_x\geq 2\) then \(f^{n_x-1}(f(x))=a\) and so \[\phi(f(x))=f_{C_k}^{-n_{f(x)}}(0)=f_{C_k}^{-n_x+1}(0)=f_{C_k}(f_{C_k}^{-n_x}(0))=f_{C_k}(\phi(x)).\]
If \(n_x=1\) then \[\phi(f(x))=\phi(a)=0=f_{C_k}(f_{C_k}^{-1}(0))=f_{C_k}(\phi(x)).\]
If \(x=a\) then \(n_{f(x)}=k-1\), and so \[\phi(f(x))=f_{C_k}^{-k+1}(0)=f_{C_k}^k(f_{C_k}^{-k+1}(0))=f_{C_k}(0)=f_{C_k}(\phi(x)).\qedhere\]
\end{proof}
\end{lemma}

\begin{lemma} \label{LineHomom}
Let \(A\) be a monounary algebra. Suppose \(a\in A\) such that \(f^{-n}(a)=\varnothing\) for some \(n\in\mathbb{N}\). Define \(\lambda_a:A\to L_{n+1}\) by: 
\begin{equation*}
    \lambda_a(x)=
    \begin{cases}
        m+1 & \text{if }x\in f^{-m}(a),\\
        0 & \text{else}.
        \end{cases}
\end{equation*}
Then
\begin{enumerate}[label=\roman*)]
    \item \(\lambda_a\) is a homomorphism.
    \item \(\lambda_a(a)=\lambda_a(b)\) if and only if \(a=b\).
\end{enumerate}
\begin{proof}
By Lemma \ref{PreimageLemmas}, \(a\) is not in a cycle, thus the \(f^{-i}(a)\) are  disjoint for all \(i<n\), so \(\lambda_a\) is a well-defined function. If \(x\in f^{-k}(a)\) for some \(k\geq1\) then by Lemma \ref{PreimageLemmas}, we have \(\lambda_a(f(x))=\lambda_a(x)-1=f_{L_{n+1}}(\lambda_a(x))\). If \(x=a\) then \(f(x)=f(a)\not\in f^{-k}(a)\) for any \(k\in\mathbb{N}\), so \(\lambda_a(f(x))=0=\lambda_a(x)-1=f_{L_{n+1}}(\lambda_a(x))\). Finally, if \(x\not\in f^{-k}(a)\) for all \(k\in\mathbb{N}\) then by Lemma \ref{PreimageLemmas}, \(f(x)\not\in f^{-k}(a)\) for all \(k\in\mathbb{N}\). Hence \(\lambda_a(f(x))=0=f_{L_{n+1}}(\lambda_a(x))\). Thus \(\lambda_a\) is a homomorphism.
Furthermore, \(\lambda_a^{-1}(\lambda_a(a))=\lambda_a^{-1}(1)=f^{-0}(a)=\{a\}\).
\end{proof}
\end{lemma}

\begin{notation}
This type of homomorphism will make repeated appearances throughout the paper, as it allows us to separate \(a\) from all other elements of \(A\) in a finite homomorphic image. As such we will reserve the notation \(\lambda_a\) specifically for these homomorphisms.
\end{notation}

\section{A Graphical Characterisation of Residual Finiteness} \label{sec:RFCrit}

This section provides a criterion for residual finiteness of monounary algebras. The proof of the criterion uses results from previous sections, in particular Lemmas \ref{SeparableBackpaths} and \ref{UniqueCyclesLemma}. But first, we briefly discuss how connectedness can affect residual finiteness.

\begin{lemma} \label{ConnectedComponentsRisdualFiniteness}
Let \((A,f)\) be a monounary algebra, with connected components \(\{K_i:i\in\mathcal{I}\}\) (where \(\mathcal{I}\) is an arbitrary index set). Then \(A\) is residually finite if and only if \(K_i\) is residually finite for all \(i\in\mathcal{I}\).
\begin{proof}
First assume each of the \(K_i\) is residually finite. Let \(x\),\(y\in A\), then there exist \(i_x\),\(i_y\in\mathcal{I}\) such that \(x\in K_{i_x}\) and \(y\in K_{i_y}\). Then if \(i_x\neq i_y\) we can construct a homomorphism \(\phi:A\to T_2\) as follows:
\begin{equation*}
    \phi(a)=
    \begin{cases}
        0 & \text{if }a\in K_{i_x},\\
        1 & \text{otherwise}.
    \end{cases}
\end{equation*}
It is easily verifiable that \(\phi\) is a homomorphism, and \(\phi(x)\neq\phi(y)\).
If \(i_x=i_y\) then since \(K_{i_x}\) is residually finite, there exists a homomorphism \(\phi:K_{i_x}\to F\), where \(F\) is finite, such that \(\phi(x)\neq\phi(y)\). Let \(F'=F\sqcup T_1\) be the disjoint union of \(F\) with \(T_1\). Then we extend \(\phi\) to a homomorphism \(\Phi:A\to F'\) by letting \(\Phi(A\setminus K_{i_x})= 0\) where \(0\) is the single point in \(T_1\).
Thus \(A\) is residually finite.

For the converse, simply note that residual finiteness is preserved under taking subalgebras.
\end{proof}
\end{lemma}

\begin{lemma} \label{inseperableBEpoints}
Let \((A,f)\) be a monounary algebra, and \((B,g)\) a finite monounary algebra such that there is a homomorphism \(\phi:A\to B\). If \(x,y\in A\) are distinct backwards eternal elements such that \(f(x)=f(y)\), then \(\phi(x)=\phi(y)\).

\begin{proof}
As \(x,y\) are in the same connected component of \(A\), we can assume without loss of generality that \(A\) is connected. We may also assume without loss of generality that \(\phi\) is surjective, and thus that \(B\) is connected. 

Since \(f^{-n}(x)\neq\varnothing\) for all \(n\in\mathbb{N}\), then there exists \(x_n\in A\) such that \(f^n(x_n)=x\). Applying \(\phi\), we get \(\phi(f^n(x_n))=g^n(\phi(x_n))=\phi(x)\). Thus \(g^{-n}(\phi(x))\neq\varnothing\) for all \(n\in\mathbb{N}\). However, since \(B\) is finite, there must be elements which lie in the intersection of two preimage sets of \(\phi(x)\), and so it follows from Lemma \ref{PreimageLemmas} \textit{iv)} that \(\phi(x)\) is in the unique cycle in \(B\).
However, we can apply the same logic to \(y\) to see that \(\phi(y)\) is also in the cycle. 
But \(g(\phi(x))=\phi(f(x))=\phi(f(y))=g(\phi(y))\), and since \(g\) is a bijection when restricted to the cycle, it follows that \(\phi(x)=\phi(y)\).

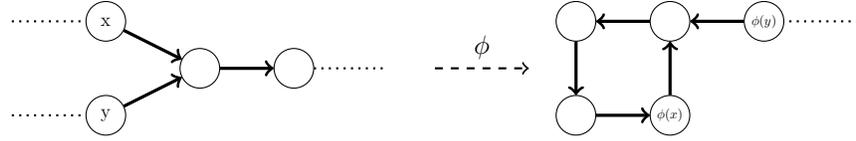
\begin{figure}
    \centering
    \begin{tikzpicture}[scale=1.25]
        \foreach \x/\place/\n in {{1/(-3,-0.5)/y},{2/(-3,0.5)/x},{3/(-2,0)/},{4/(-1,0)/}}
            \node[V,minimum size=30] (z\x) at \place {\Large{\n}};
        \draw[E] (z1) -- (z3);
        \draw[E] (z3) -- (z4);
        \draw[E] (z2) -- (z3);
        \draw[dotted, thick] (-4,0.5) -- (z2);
        \draw[dotted, thick] (z4) -- (0,0);
        \draw[dotted, thick] (-4,-0.5) -- (z1);
        \foreach \x/\place in {{1/(2,-0.5)},{2/(3,-0.5)},{3/(3,0.5)},{4/(2,0.5)},{5/(4,0.5)}}
            \node[V,minimum size=30] (c\x) at \place {};
        \node[scale=0.5] (lpy) at (4,0.5) {\(\phi(y)\)};
        \node[scale=0.5] (lpx) at (3,-0.5) {\(\phi(x)\)};
        \draw[E] (c1) -- (c2);
        \draw[E] (c2) -- (c3);
        \draw[E] (c3) -- (c4);
        \draw[E] (c4) -- (c1);
        \draw[E] (c5) -- (c3);
        \draw[dotted, thick] (5,0.5) -- (c5);
        \draw[dashed,->, thick] (0.5,0) -- (1.5,0);
        \node[above] (phi) at (1,0) {\(\phi\)};
    \end{tikzpicture}
    \caption{Showing the results of attempting to get a finite homomorphic image of an algebra which fails the RF criterion.}
    \label{fig:NonRFAttempt}
\end{figure}
\end{proof}
\end{lemma}

\begin{theorem} \label{MainRFThm}
Let \((A,f)\) be a monounary algebra. Then the following are equivalent: 
\begin{enumerate}[label=(\arabic*)]
    \item \(A\) is residually finite,
    \item for all distinct \(x,y\in A\) such that \(f(x)=f(y)\), there exists \(n\in\mathbb{N}\) such that either \(f^{-n}(x)=\varnothing\) or \(f^{-n}(y)=\varnothing\).
\end{enumerate} 

\begin{notation}
As we will make repeated reference to it, for simplicity we will refer to the second condition as the \textit{RF criterion}.
\end{notation}

\begin{proof}
By Lemma \ref{ConnectedComponentsRisdualFiniteness}, we may assume without loss of generality that \(A\) is connected.

\textbf{\textit{(1)}\(\Rightarrow\)\textit{(2)}} Follows immediately from Lemma \ref{inseperableBEpoints}.

\textbf{\textit{(2)}\(\Rightarrow\)\textit{(1)}}
Next we show that if the RF criterion holds, then the algebra is residually finite. We shall do this by constructing homomorphisms. Let \(x,y\in A\). Note that via Lemma \ref{SeparableBackpaths}, we have that one of the following holds 
\begin{enumerate}
    \item There exists \(n\in\mathbb{N}\) such that \(f^n(x)=y\) or \(f^n(y)=x\).
    \item \((\bigcup_{n\in\mathbb{N}}f^{-n}(x))\cap(\bigcup_{n\in\mathbb{N}}f^{-n}(y))=\varnothing\).
\end{enumerate}

We shall first deal with case 2. By Lemma \ref{SeparableBackpaths}, we may find minimal \(i,j\in\mathbb{N}\) (with \(i,j\geq1\) such that \(f^i(x)=f^j(y)\). Call this common point \(p\). It then follows that there exist \(x'=f^{i-1}(x)\neq f^{j-1}(y)=y'\) such that \(f(x')=f(y')=p\). It then follows by the RF criterion that, without loss of generality, there exists \(n'\in\mathbb{N}\) such that \(f^{-n'}(x')=\varnothing\) (and by Lemma \ref{PreimageLemmas}, \(n'>i-1\)), and so setting \(n=n'-(i-1)\) we have \(f^{-n}(x)=\varnothing\). Then we can use \(\lambda_x\) from Lemma \ref{LineHomom}, which is a homomorphism to a finite algebra, and \(\lambda_x(x)\neq\lambda_x(y)\).

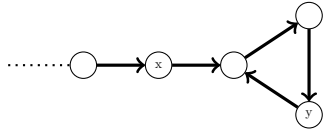
\begin{figure}
    \centering
    \begin{tikzpicture}
        \foreach \x/\place/\n in {{1/(0,0)/},{2/(1,0.66)/},{3/(1,-0.66)/y},{4/(-1,0)/x},{5/(-2,0)/}}
            \node[V,minimum size=20] (a\x) at \place {\n};
        \draw[E] (a1)--(a2);
        \draw[E] (a2)--(a3);
        \draw[E] (a3)--(a1);
        \draw[E] (a4)--(a1);
        \draw[E] (a5)--(a4);
        \draw[dotted, thick] (-3,0)--(a5);
    \end{tikzpicture}
    \caption{RF criterion fails due to $x$ and $y$.}
    \label{fig:RFCriterionFailure}
\end{figure}

For case 1, note that if there exists a finite cycle \(C_k\subset A\) and a backwards eternal element \(a\in A\setminus C_k\), then the RF criterion does not hold. Thus if there exists a cycle, then for every point \(a\in A\) which is not in the cycle, there exists \(m\in\mathbb{N}\) such that \(f^{-m}(a)=\varnothing\). We can thus further separate case 1 into subcases:
\begin{enumerate}[label=(\alph*)]
    \item both \(x\) and \(y\) are in the cycle \(C_k\);
    \item there exists \(m\in\mathbb{N}\) such that \(f^{-m}(x)=\varnothing\) or \(f^{-m}(y)=\varnothing\);
    \item there are no cycles in \(A\).
\end{enumerate}

For subcase 1(a), note that as the algebra is connected and contains a cycle we can use the homomorphism \(\phi\) defined in Lemma \ref{cyclehomom}. Then note that this homomorphism separates all elements of the cycle from each other, and so \(\phi(x)\neq\phi(y)\).

For subcase 1(b), we can again use \(\lambda_x\) or \(\lambda_y\).

For subcase 1(c), we note that if we are not also in subcase 1(b), we have that \(y\) is backwards eternal and \(f^i(y)\neq f^j(y)\) for all distinct \(i, j\in\mathbb{N}\). Let \(Z=\{f^i(y):i\in\mathbb{N}\}\cup(\bigcup_{n\in\mathbb{N}}f^{-n}(y))\). Let \(z\in A\setminus Z\), then by Lemma \ref{SeparableBackpaths}, there exist \(i\),\(j\in\mathbb{N}\) such that \(f^i(z)=f^j(y)\). Thus for each \(a\in A\setminus Z\) we associate the value \(p(z)=j-i\in\mathbb{Z}\). Then let \(P_k=\{z\in A\setminus Z:p(z)=k\}\). We thus have that \(A=Z\sqcup(\bigsqcup_{i\in\mathbb{Z}}P_i)\). We then define a function \(\theta:A\to \mathbb{Z}\) as follows:
\begin{equation*}
    \theta(a)=
    \begin{cases}
        -k & \text{if }a\in f^{-k}(y)\cup P_{-k},\\
        0 & \text{if }a\in \{y\}\cup P_0,\\
        k & \text{if }a\in \{f^k(y)\}\cup P_k.
    \end{cases}
\end{equation*}

We show that this is a homomorphism to the monounary algebra \(\mathbb{Z}\) and that \(\theta(x)\neq\theta(y)\).
If \(a\in f^{-k}(y)\) for some \(k\geq 1\) then \(f(a)\in f^{-k+1}(y)\), and so \(\theta(f(a))=-k+1=\theta(a)+1\).
If \(a=y\) then \(\theta(f(a))=1=0+1=\theta(a)+1\).
If \(a=f^k(y)\) for some \(k\geq 1\) then \(\theta(f(a))=\theta(f^{k+1}(y))=k+1=\theta(a)+1\).
If \(a\in P_k\) for some \(k\in\mathbb{Z}\) and \(f(a)=f^{k+1}(y)\) then \(\theta(f(a))=k+1=\theta(a)+1\).
If \(a\in P_k\) for some \(k\in\mathbb{Z}\) and \(f(a)\neq f^{k+1}(y)\) then \(f(a)\in P_{k+1}\) so \(\theta(f(a))=k+1=\theta(a)+1\).
Thus \(\theta\) is a homomorphism. And since we have that \(f^n(x)=y\) or \(f^n(y)=x\) for some \(n\in\mathbb{N}\), it is clear that \(\theta(x)\neq\theta(y)\).

Then using Lemma \ref{ZRF}, we can construct a homomorphism \(\sigma:\mathbb{Z}\to C_m\) for some \(m\in\mathbb{N}\) such that \(\sigma(\theta(x))\neq\sigma(\theta(y))\). Thus the composition \(\sigma\theta\) is a homomorphism into a finite algebra which separates \(x\) and \(y\) as required.
\end{proof}
\end{theorem}

\begin{figure}[h]
    \centering
    \begin{tikzpicture} [scale=1.25]
        \foreach \x/\place in {{-2/(-2,0)},{-1/(-1,0)},{0/(0,0)},{-2i/(-2,-1)},{-1i/(-1,-1)},{2/(2,0)},{1/(1,0)}}
            \node[V] (b\x) at \place {};
        \foreach \o/\i in {{-2/-1},{-1/0},{-2i/-1i},{-1i/0},{0/1},{1/2}}
            \draw[E] (b\o) -- (b\i);
        \draw [dotted, thick] (-3,0) -- (b-2);
        \draw [dotted, thick] (-3,-1) -- (b-2i);
        \draw [dotted, thick] (b2) -- (3,0);
    \end{tikzpicture}
    \caption{An example of a monounary algebra which is not residually finite.}
    \label{fig:NonRFExample}
\end{figure}
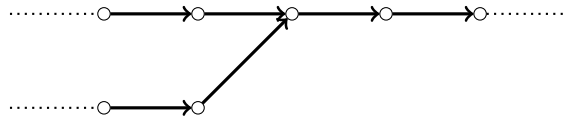

If we were to rephrase the RF criterion in graphical terms, it would be as follows:``For all distinct vertices \(x,y\), if there exists a vertex \(z\) such that \((x,z),(y,z)\) are edges, then at least one of \(x\) or \(y\) has a finite bound on the length of walks that terminate there." We can also phrase it in terms of backwards eternality, and it would become: ``For all \(a\in A\), at most one element of \(f^{-1}(a)\) is backwards eternal." 

\section{Direct Products} \label{sec:Prod}

Our primary goal in this section is to use our criterion from Section~\ref{sec:RFCrit} to obtain necessary and sufficient conditions on components of a direct product of monounary algebras for the direct product itself to be residually finite.

We first take note of a class of monounary algebras that seems to exhibit residual finiteness that is, in some sense, more powerful than usual. We shall give this particular type of monounary algebra a name, as they turn out to have some very strong separation properties, and behave in fundamentally different ways with respect to direct products.

\begin{definition}
A monounary algebra \(A\) is called \textit{backwards-bounded} if for all \(a\in A\) there exists an \(n\in\mathbb{N}\) such that \(f^{-n}(a)=\varnothing\). Equivalently, it is a mmonounary algebra which contains no backwards eternal elements.
\end{definition}

Note that it is possible for \(f^{-m}(a)\) to be an infinite set, as there could be infinitely many paths ending at \(a\), but with finite maximum length.

Rephrased in graphical terms, this becomes: ``A monounary algebra \(A\) is backwards-bounded if for all vertices \(x\) there is a finite bound on the length of walks which terminate at \(x\)."

By definition such structures satisfy the RF criterion, and so;
\begin{lemma}
Backwards-bounded monounary algebras are residually finite.
\end{lemma}

In addition, for a direct product \((A\times B,f)\) of monounary algebras \((A,f_1),(B,f_2)\), we have \(f^n((a,b))=(f_1^n(a),f_2^n(b))\). From this we get the following lemma, which shows that the finite number of non-empty preimages is particularly powerful. From here on, we will drop the double brackets and write \(f(a,b)\) instead of \(f((a,b))\), and whenever we have a direct product, we will use \(f_i\) to refer to the operation on the \(i\)th component and \(f\) to refer to the operation on the product unless otherwise specified.

\begin{lemma} \label{ProductBackpathLemma}
For a direct product \(A\times B\) of monounary algebras \(A\) and \(B\), we have \(f^{-n}(x,y)=f_1^{-n}(x)\times f_2^{-n}(y)\).
\begin{proof}
\begin{align*}
    (a,b)\in f^{-n}(x,y) & \Leftrightarrow f^n(a,b)=(x,y),\\ &\Leftrightarrow f_1^n(a)=x, f_2^n(b)=y\\ &\Leftrightarrow a\in f_1^{-n}(x) , b\in f_2^{-n}(y). 
\end{align*}
\end{proof}

\end{lemma}

This yields the following two propositions about residual finiteness of certain direct products, which essentially show that backwards-boundedness forces direct products to be residually finite, and that it is the only class of monounary algebras that do so.

\begin{proposition} \label{StrongRFProp}
Let \(A\) be a backwards-bounded monounary algebra. Then for any monounary algebra \(B\) we have that \(A\times B\) is backwards bounded, and hence residually finite.
\begin{proof}
Let \((a,b)\in A\times B\). Then there exists \(k_a\in\mathbb{N}\) such that \(f_1^{-k_a}(a)=\varnothing\). It then follows by Lemma \ref{ProductBackpathLemma} that \(f^{-k_a}(a,b)=\varnothing\). Thus we have that for every point \(x\in A\times B\), there exists an \(n\in\mathbb{N}\) such that \(f^{-n}(x)=\varnothing\), and so \(A\times B\) is backwards-bounded. Thus \(A\times B\) is residually finite.
\end{proof}
\end{proposition}

\begin{proposition} \label{backwards-boundedRFprop}
Let \(A\) be a non-residually finite monounary algebra, and \(B\) another monounary algebra. If \(A\times B\) is residually finite then \(B\) is backwards-bounded.
\begin{proof}
We show the contrapositive: assume \(B\) is not backwards-bounded, so there exists \(b\in B\) which is backwards eternal.

As \(A\) is not residually finite, the RF criterion does not hold, so there exist two distinct backwards eternal points \(x,y\in A\). Consider the points \((x,b),(y,b)\in A\times B\). These are distinct points and \(f(x,b)=(f_1(x),f_2(b))=(f_1(y),f_2(b))=f(y,b)\). Then since \(f_1^{-n}(x),f_1^{-n}(y),f_2^{-n}(b)\neq\varnothing\) for all \(n\in\mathbb{N}\), it follows by Lemma \ref{ProductBackpathLemma} that \(f^{-n}(x,b),f^{-n}(y,b)\neq\varnothing\) for all \(n\in\mathbb{N}\). Thus by Theorem \ref{MainRFThm}, \(A\times B\) is not residually finite.
\end{proof}
\end{proposition}

These propositions combine to give us the following theorem, which determines the residual finiteness of a direct product of monounary algebras from the properties of the components. In particular it shows there is only one way to get a residually finite product without both components being residually finite. 

\begin{theorem} \label{RFproducttheorem}
A direct product of monounary algebras, \(A\times B\), is residually finite if and only if one of the following holds:
\begin{enumerate}
    \item Both \(A\) and \(B\) are residually finite.
    \item \(A\) is backwards-bounded.
    \item \(B\) is backwards-bounded.
\end{enumerate}
\begin{proof}
The reverse implication follows trivially for 1, and from Proposition \ref{StrongRFProp} for 2 and 3.
For the forward implication, let \(A\times B\) be residually finite. Then if 1 is not true, then at least one of \(A\) and \(B\) is not residually finite, so by Proposition \ref{backwards-boundedRFprop}, the other is backwards-bounded, giving 2 or 3.
\end{proof}
\end{theorem}

We can easily extend this result to arbitrary products.

\begin{theorem}
For an arbitrary index set \(\mathcal{I}\), a direct product of monounary algebras \(\prod_{i\in I}X_i\) is residually finite if and only if one of the following holds:
\begin{enumerate}
    \item \(X_i\) is residually finite for all \(i\),
    \item There exists an \(i\) such that \(X_i\) is backwards-bounded.
\end{enumerate}
\begin{proof}
First note that it can be seen shortly from Lemma \ref{ProductBackpathLemma} that in fact a direct product is backwards bounded if and only if at least one of its component is backwards-bounded.

For the reverse direction, note we can show 1 implies the product is residually finite in exactly the same way as in the case for two factors. For 2, as at least one component is backwards-bounded, the product is backwards-bounded (and thus residually finite).
For the forward implication, note that if at least one of the factors \(X_i\) is not residually finite and the rest are not backwards-bounded, then the remaining product \(\prod_{j\in\mathcal{I}\setminus\{i\}}X_j\) is not backwards bounded, and so the whole product is not residually finite. Thus if the product is residually finite, and at least one of the factors is not residually finite, then at least one factor must be backwards-bounded.
\end{proof}
\end{theorem}

As a brief aside, we consider \textit{subdirect products}, subalgebras of the direct product such that the projection maps are surjective. We show that the equivalence we obtained for direct products does not hold, by constructing an explicit example of a residually finite subdirect product of two monounary algebras, neither of which are residually finite.

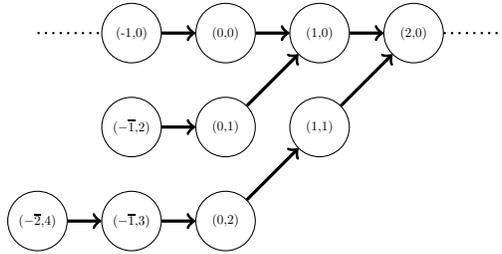
\begin{figure}
    \centering
    \begin{tikzpicture}[scale=1.25]
        \foreach \x/\place/\n in {{-10/(-1,0)/(-1,0)},{00/(0,0)/(0,0)},{10/(1,0)/(1,0)},{20/(2,0)/(2,0)},{01/(0,-1)/(0,1)},{i2/(-1,-1)/(\(-\overline{1}\),2)},{11/(1,-1)/(1,1)},{02/(0,-2)/(0,2)},{1i3/(-1,-2)/(\(-\overline{1}\),3)},{2i4/(-2,-2)/(\(-\overline{2}\),4)}}
            \node[V,minimum size=45] (s\x) at \place {\n};
        \foreach \o/\i in {{-10/00},{00/10},{10/20},{i2/01},{01/10},{2i4/1i3},{1i3/02},{02/11},{11/20}}
            \draw[E] (s\o) -- (s\i);
            \draw[dotted, thick] (-2,0) -- (s-10);
            \draw[dotted, thick] (s20) -- (3,0);
    \end{tikzpicture}
    \caption{A residually finite subdirect product of not residually finite monounary algebras.}
    \label{fig:RFsubdirectWxample}
\end{figure}

\begin{example} \label{Subdirectproductexample}
We define a monounary algebra \(A\) on the set \(\mathbb{Z}\cup\{-\overline{n}:n\in\mathbb{N}^+\}\) by
\begin{equation*}
f(x) =
    \begin{cases}
        0 & \text{if }x=-\overline{1},\\
        -\overline{(n-1)} &\text{if }x=-\overline{n}, n>1,\\
        x+1 &\text{otherwise}.
    \end{cases}  
\end{equation*}
Note: What we are doing is attaching another disjoint copy of the negatives to the integers, and this is the same monounary algebra that is depicted in Figure \ref{fig:NonRFExample}.

We then consider the direct product of this algebra with \(N=(\mathbb{N},\max\{x-1,0\})\). Both components fail the RF criterion (with the points \(-1\) and \(-\overline{1}\) for \(A\) and \(0\) and \(1\) for \(B\)) so by Theorem \ref{RFproducttheorem}, \(A\times N\) is not residually finite. However, we can construct a subdirect product which is residually finite. Consider the subalgebra of the direct product \(A\times N\) generated by the set \(G=\{(n,0):n\in\mathbb{Z}\}\cup\{(-\overline{m},2m):m\in\mathbb{N}^+\}\) (this is depicted in figure \ref{fig:RFsubdirectWxample}). This is clearly a subdirect product, and satisfies our criterion from Theorem \ref{MainRFThm}, as each element \(x\) other than those in \(\mathbb{Z}\times\{0\}\) has the property that \(f^{-n}(x)=\varnothing\) for a large enough \(n\in\mathbb{N}\).
\end{example}

It is worth noting however, that as backwards-boundedness is preserved under subalgebras, if either component is backwards-bounded, then any subdirect product is also backwards-bounded, and hence residually finite.

\section{Further Separability Properties}

In this section we discuss three notions related to residual finiteness: weak and strong subalgebra separability, and complete separability. While these properties have been studied for some time, such as in \cite{mal1983homomorphisms} and \cite{golubov1970finite}, the names used for them have not been consistent. The names we use are from \cite{miller2020separability}, as these are designed to be more descriptive of the property. We provide characterisations for these properties, and show how they interact with direct products in a similar fashion to how we dealt with residual finiteness.

\begin{definition}
A monounary algebra \(A\) is \textit{strongly (weakly) subalgebra separable} if for any \(a\in A\) and any (any finitely generated) subalgebra \(B\leq A\) such that \(a\not\in B\), there exists a finite monounary algebra \(F\) and a homomorphism \(\phi:A\to F\) such that \(\phi(a)\not\in\phi(B)\).
\end{definition}

To deal with these conditions we will introduce the notion of bi-eternal monounary algebras. We show in Theorem \ref{SSEquivalence} that these algebras are the only ones which distinguish residual finiteness from strong and weak subalgebra separability. 

\begin{definition}
A monounary algebra \(A\) is \textit{bi-eternal} if there exists \(a\in A\) such that the following two conditions hold:
\begin{enumerate}
    \item \(f^i(a)=f^j(a)\) if and only if \(i=j\),
    \item \(a\) is backwards eternal.
\end{enumerate}
\end{definition}

\begin{remark}
From here on we shall refer to the first criterion from the above definition as the forward eternality.
\end{remark}

Rephrased in graphical terms, the forward eternality condition becomes: ``The unique infinite walk starting at \(a\) is an infinite path."

\begin{example}
The monounary algebra \(\mathbb{Z}\) is bi-eternal, as any point satisfies both eternality conditions.
\end{example}

Example \ref{nonZbieternal} is an example of a bi-eternal monounary algebra which does not contain \(\mathbb{Z}\) as a subalgebra.

It is useful to note that the forward eternality condition corresponds to not having a cycle (in the corresponding connected component) and the backwards eternality condition means we do not have backwards-boundedness. We can thus split connected monounary algebras into three distinct classes.

\begin{lemma} \label{BEsplit}
Let \(A\) be a connected monounary algebra. Then exactly one of the following is true:
\begin{enumerate}
    \item \(A\) contains a cycle,
    \item \(A\) is bi-eternal,
    \item \(A\) is backwards bounded.
\end{enumerate}
\begin{proof}Assume \(A\) is not bi-eternal. If there is a point for which the forward eternality condition fails, then there exists a cycle in \(A\), and since \(A\) is connected, every point therefore fails the forward eternality condition. Otherwise, the forward eternality condition holds for every point, thus no point is backwards eternal, and so \(A\) is backwards bounded.
\end{proof}
\end{lemma}

\begin{remark}
Note that, in much the same way as residual finiteness, weak and strong subalgebra separability hold if and only if they hold for every connected component. 
\end{remark}

We now show that strong and weak subalgebra separability are equivalent and that bi-eternality is the only thing which separates them from residual finiteness.

\begin{theorem} \label{SSEquivalence}
For a monounary algebra \(A\), the following are equivalent:
\begin{enumerate}[label=(\arabic*)]
    \item \(A\) is strongly subalgebra separable,
    \item \(A\) is weakly subalgebra separable,
    \item \(A\) is residually finite and not bi-eternal,
    \item For all \(x\in A\) either \(x\) is contained in a cycle or there exists \(n\in\mathbb{N}\) such that \(f^{-n}(x)=\varnothing\).
\end{enumerate}
\begin{proof}
Since bi-eternality is determined by the existence of a point with certain properties, a monounary algebra is bi-eternal if and only if at least one of its connected components is bi-eternal, and so we can assume without loss of generality that \(A\) is connected.

\textbf{\textit{(1)}\(\Rightarrow\)\textit{(2)}}. 
This follows from the definitions.

\textbf{\textit{(2)}\(\Rightarrow\)\textit{(3)}}
We show the contrapositive.

If \(A\) is not residually finite, then we can take the two points which fail the RF criterion, \(x\) and \(y\). Since they are distinct points with the same image under the unary operation, at most one of them can be in a cycle. Thus we must have at least one of \(x\not\in\langle y\rangle\) or \(y\not\in\langle x\rangle\). But as these two points are backwards eternal, by Lemma \ref{inseperableBEpoints} they cannot be mapped to distinct points in a finite algebra , and so we must have that \(\phi(x)\in\langle \phi(y)\rangle\) and \(\phi(y)\in\langle \phi(x)\rangle\) for any \(\phi\) a homomorphism from \(A\) to a finite monounary algebra. Hence \(A\) is not weakly subalgebra separable.

If \(A\) is bi-eternal, then consider the element \(a\) for which the eternality conditions hold. The forward eternality condition shows that \(a\) is not in a cycle, and so \(a\not\in\langle f(a)\rangle\). Now if \(\phi\) is a homomorphism from \(A\) to a finite monounary algebra \((F,g)\), then \(\phi(\langle f(a)\rangle)\) is a non-empty subalgebra of \(F\) and so contains the cycle of \(F\). But as \(a\) is backwards eternal, for all \(n\in\mathbb{N}\) there exists \(a_n\in A\) such that \(g^n(\phi(a_n))=\phi(a)\). As \(F\) is finite, some of these must be the same, which forces \(\phi(a)\) to be in the cycle of \(F\), so \(\phi(a)\in\phi(\langle f(a)\rangle)\). Thus \(A\) is not weakly subalgebra separable.

\textbf{\textit{(3)}\(\Rightarrow\)\textit{(1)}}
Since \(A\) is not bi-eternal, we can use Lemma \ref{BEsplit} to split into cases:
\begin{enumerate}
    \item \(A\) is backwards-bounded,
    \item \(A\) contains a cycle (and is residually finite).
\end{enumerate}

In case 1, for any \(a\in A\) and \(B\subset A\) such that \(a\not\in B\), we can use \(\lambda_a\) as defined in Lemma \ref{LineHomom}, as in this case \(\phi(a)\not\in\phi(A\setminus\{a\})\supseteq \phi(B)\).

In case 2, since the cycle is contained in every non-empty subalgebra, if \(a\in A\) is in the cycle of \(A\) then \(B=\varnothing\), and we are done. If \(a\in A\) is not in the cycle, then \(f^{-n}(a)=\varnothing\) for some \(n\in\mathbb{N}\). Thus we can once again use \(\lambda_a\).

Thus in either case \(A\) is strong subalgebra separable. 

\textbf{\textit{(3)}\(\Rightarrow\)\textit{(4)}}
Using Lemma \ref{BEsplit} we can conclude that \(A\) is either backwards-bounded (in which case we are done), or contains a cycle and is residually finite. But since every element \(a\) of the cycle has \(f^{-n}(a)\neq\varnothing\) for all \(n\), in order to be residually finite, we must have that for every element \(x\) outside the cycle \(f^{-n}(x)=\varnothing\) for some \(n\in\mathbb{N}\).

\textbf{\textit{(4)}\(\Rightarrow\)\textit{(3)}}
Now assume for every \(x\in A\) either \(x\) is contained in a cycle or there exists \(n\in\mathbb{N}\) such that \(f^{-n}(x)=\varnothing\). Thus the only backwards eternal elements are the cycle elements (if the cycle exists) and so the RF criterion is satisfied, and it is clear to see by the definition that \(A\) is not bi-eternal, as the forward eternality condition is not satisfied by any point.
\end{proof}
\end{theorem}

Rephrasing the fourth condition in graphical term, it becomes: ``For any vertex \(x\) which is not in a cycle, then there is a finite bound on the length of walks terminating at \(x\)." In terms of backwards eternality, it is: ``Every backwards eternal element is contained in a cycle."

For the rest of this document we shall refer to such algebras as subalgebra separable for simplicity.

\begin{lemma} \label{backwardsboundeddominateszl}
Let \(A\) be a bi-eternal monounary algebra. Then for \(B\) a monounary algebra , \(A\times B\) is bi-eternal if and only if \(B\) is not backwards-bounded.
\begin{proof}
The forward implication follows immediately from Proposition \ref{StrongRFProp}.

For the converse, note that if \(B\) is not backwards-bounded then it is either bi-eternal or contains a cycle. If \(x\in A\) satisfies the eternality conditions, then in both cases we identify a corresponding point in the direct product which satisfies the eternality condition. If \(B\) contains a cycle, let \(b\) be a point in the cycle. Then \((x,b)\) satisfies the eternality conditions: the first is inherited from \(x\), and the second follows from Lemma \ref{ProductBackpathLemma}. If \(B\) is instead bi-eternal, then we can take a point \(y\in B\) which satisfies the eternality conditions, and \((x,y)\) also satisfies the eternality conditions by the same argument as for \((x,b)\). 
\end{proof}
\end{lemma}

\begin{lemma} \label{zlproductlemma}
For monounary algebras \(A\) and \(B\), if \(A\times B\) is bi-eternal, then at least one of \(A\) and \(B\) is bi-eternal.
\begin{proof} 
Since a monounary algebra is bi-eternal if at least one connected component is bi-eternal, we may assume without loss of generality that \(A\) and \(B\) are connected.
Note that \(A\times B\) is backwards-bounded if and only if at least one of its components is backwards-bounded. Thus neither \(A\) nor \(B\) is backwards bounded and so by Lemma \ref{BEsplit} each must be either bi-eternal or contain a cycle. Let us assume for a contradiction, that both contain a cycle. Then for every \((a,b)\in A\times B\), there exists \(n,m\in\mathbb{N}\) such that \(f_1^n(a)\) is in the cycle of \(A\) and \(f_2^m(b)\) is in the cycle of \(B\). Then, since a pair with both coordinates in the corresponding cycle of the component is in a cycle in the product, \(f^{n+m}(a,b)\) is in a cycle, and so \((a,b)\) does not satisfy the eternality condition. Since \((a,b)\) was arbitrary, \(A\times B\) is not bi-eternal, a contradiction. 
\end{proof}
\end{lemma}

Using these results, we obtain a result that mirrors Theorem \ref{RFproducttheorem}. In particular, the conditions for a direct product of monounary algebras to be subalgebra separable are the same as those for being residually finite (replacing residually finite with subalgebra separable).

\begin{theorem} \label{SSproducttheorem}
Let \(A\) and \(B\) be monounary algebras. Then \(A\times B\) is subalgebra separable if and only if one of the following is true.
\begin{enumerate}
    \item \(A\) and \(B\) are subalgebra separable,
    \item \(A\) is backwards-bounded,
    \item \(B\) is backwards-bounded.
\end{enumerate}

\begin{proof}
Note that \(A\times B\) is the union of the direct products of each connected component of \(A\) with each connected component of \(B\). Then since a monounary algebra is subalgebra separable if and only if its connected components are subalgebra separable (and similarly for backwards-bounded), we may assume without loss of generality that \(A\) and \(B\) are connected. 

For the converse, if one component is backwards-bounded, then by Proposition \ref{StrongRFProp}, \(A\times B\) is backwards-bounded, and thus subalgebra separable by Theorem \ref{SSEquivalence}, item (4). If both \(A\) and \(B\) are subalgebra separable, then they are both residually finite and not bi-eternal, and so by Theorem \ref{RFproducttheorem} and Lemma \ref{zlproductlemma}, the product is both residually finite and not bi-eternal, and thus is subalgebra separable.

For the forward implication, we show the contrapositive. Thus we assume one of \(A\) or \(B\) (say \(A\)) is not subalgebra separable (and so by Theorem \ref{SSEquivalence}, either is not residually finite or is bi-eternal), and that neither is backwards-bounded. In particular, as neither are backwards-bounded, Theorem \ref{RFproducttheorem} tells us that if \(A\) was not residually finite, the direct product is not residually finite, and thus not subalgebra separable by Theorem \ref{SSEquivalence}. If \(A\) was instead bi-eternal then, by Lemma \ref{backwardsboundeddominateszl} the direct product is bi-eternal and thus not subalgebra separable.
\end{proof}
\end{theorem}

Having dealt with subalgebra separability, we now move on to the concept of complete separability.

\begin{definition}
A monounary algebra is \textit{completely separable} if for every \(a\in A\) there exists a finite monounary algebra \(F\) and  homomorphism \(\phi:A\to F\) such that \(\phi(a)\not\in\phi(A\setminus{a})\).
\end{definition}

\begin{lemma}
Let \(A\) be a monounary algebra. If \(A\) is completely separable, then it is also subalgebra separable (and thus residually finite).
\begin{proof}
This is immediate from the definitions.
\end{proof}
\end{lemma}

In fact, the converse of this lemma only fails in a specific scenario. The obstacle is cycles, which need not be considered for subalgebra separability as they are contained in every non-empty subalgebra. As such we will use the following definition to deal with complete separability.

\begin{definition}
Let \((A,f)\) be a monounary algebra. Then we define \[B_n(a)=f^{-n}(a)\setminus\bigcup_{i=0}^{n-1}f^{-i}(a)\]
\end{definition}

\begin{theorem} \label{CStheorem}
A monounary algebra \((A,f)\) is completely separable if and only if for all \(a\in A\) there exists \(n\in\mathbb{N}\) such that \(B_n(a)=\varnothing\).

\begin{proof}
As per usual, we may assume without loss of generality that \(A\) is connected.
First assume that \(A\) is completely separable. Let \(x\in A\) with corresponding finite algebra \((F,g)\) and homomorphism \(\phi:A\to F\). Consider \(\phi(x)\). If \(\phi(x)\) is not in the cycle of \(F\), then as \(F\) is finite there exists \(n\in\mathbb{N}\) such that \(g^{-n}(\phi(x))=\varnothing\). But since \(\phi\) is a homomorphism, if \(y\in f^{-n}(x)\) then \(\phi(y)\in g^{-n}(\phi(x))\). Thus \(f^{-n}(x)=\varnothing\) (and so \(B_n(x)=\varnothing\)). Now if \(\phi(x)\) is in the cycle of \(F\) then, denoting the cycle length by \(k\), we have \(\phi(f^k(x))=g^k(\phi(x))=\phi(x)\), and so by complete separability we have that \(f^k(x)=x\), and so \(x\) is in the cycle of \(A\). As \(F\) is finite, there exists \(n\in\mathbb{N}\) such that \(B_n(\phi(x))=\varnothing\). Thus if \(B_n(x)\neq\varnothing\), then \(\phi(B_n(x))\) is in the cycle. Let \(m\) be the least multiple of the cycle length greater than \(n\). Then if \(B_m(x)\neq\varnothing\), we have \(\phi(B_m(x))=\phi(x)\) a contradiction. Thus \(B_m(x)=\varnothing\).

Now assume that for all \(a\in A\) there exists \(n\in\mathbb{N}\) such that \(B_n(a)=\varnothing\).  For \(x\in A\), if \(x\) is not in a cycle, then \(B_n(x)=f^{-n}(x)=\varnothing\) so we can use the standard homomorphism \(\lambda_x\) to separate \(x\). Now let \(x\in A\) be in the cycle, and say the cycle has length \(k\). As \(A\) is connected, the cycle is the minimal non-empty subalgebra, and so every element of \(A\) is in \(B_m(x)\) for precisely one \(m\). Let \(N\in\mathbb{N}\) be the minimal value such that \(B_N(x)=\varnothing\), then as \(f(B_n(a))\subseteq B_{n-1}(a)\), it follows that \(B_M(x)=\varnothing\) for all \(M\geq N\). Thus we construct a homomorphism to an algebra \(F=(\{0,\dots,N-1\},g)\) where \(g\) is given by:

\begin{equation*}
    g(y) =
    \begin{cases}
        k &\text{if }y=0,\\
        y-1 &\text{otherwise}.
    \end{cases}
\end{equation*}

and the homomorphism \(\phi:A\to F\) is given by \(\phi(B_n(x))=n\). This is a homomorphism as \(f(B_n(x))\subseteq B_{n-1}(x)\) for \(n\geq1\) and \(f(x)\in B_k(x)\). Additionally, \(\phi^{-1}(0)=\{x\}\) and so \(x\) is completely separated by \(\phi\) as required.
\end{proof}
\end{theorem}

Rephrasing the condition in graphical terms, it becomes: ``For all vertices \(a\), there is a finite bound on the length of paths terminating at \(a\). Note that the use of \(B_n(a)\) rather than \(f^{-n}(a)\) corresponds to considering paths instead of walks."

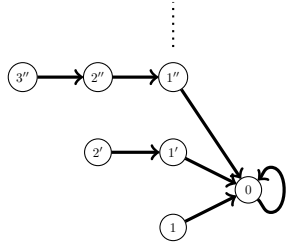
\begin{figure}
    \centering
    \begin{tikzpicture}
        \foreach \x/\place/\n in {{0/(0,0)/0},{00/(-1,-0.5)/1},{01/(-1,0.5)/$1^\prime$},{11/(-2,0.5)/2$^\prime$},{02/(-1,1.5)/$1^{\prime\prime}$},{12/(-2,1.5)/\(2^{\prime\prime}\)},{22/(-3,1.5)/\(3^{\prime\prime}\)}}
            \node[V,minimum size=20] (a\x) at \place {\n};
        \foreach \o/\i in {{22/12},{12/02},{02/0},{11/01},{01/0},{00/0}}
            \draw[E] (a\o) -- (a\i);
        \draw[dotted, thick] (-1,1.9) -- (-1,2.5);
        \draw[E] (a0) to[out=-45,in=180] (0.3,-0.3) to[out=0,in=0] (0.3,0.3) to[out=180,in=45] (a0);
    \end{tikzpicture}
    \caption{A monounary algebra which is subalgebra separable, but not completely separable.}
    \label{fig:CSexample}
\end{figure}

\begin{example}
Consider the monounary algebra formed by taking every \(n\)-line (as defined in Example \ref{exampleMUAs}) and identifying the zeros (depicted in figure \ref{fig:CSexample}). This is subalgebra separable, but not completely separable.
\end{example}

\begin{lemma}
Let \(A,B\) be monounary algebras, and \((a,b)\in A\times B\). Then \[B_n(a,b)=(B_n(a)\times f_2^{-n}(b))\cup (f_1^{-n}(a)\times B_n(b)).\]
\begin{proof}
Note that, by the definition of \(B_n\), we have that \(f^{-n}(x) = B_n(0)\sqcup\dots\sqcup B_n(x)\). Using this together with the Lemma \ref{ProductBackpathLemma} shows that the pairs in \(B_n(a,b)\) are precisely those that contain either an element of \(B_n(a)\) or an element of \(B_n(b)\).
\end{proof}
\end{lemma}

Now we get conditions for a direct product to be completely separable based on the properties of the components, and again, the conditions turn out to be the same as for residual finiteness and subalgebra separability.

\begin{theorem} \label{CSproducttheorem}
Let \(A\) and \(B\) be monounary algebras, then \(A\times B\) is completely separable if and only if one of the following is true:

\begin{enumerate}
    \item \(A\) and \(B\) are completely separable,
    \item \(A\) is backwards-bounded,
    \item \(B\) is backwards-bounded. 
\end{enumerate}

\begin{proof}
For the converse, if either \(A\) or \(B\) is backwards-bounded then by Proposition \ref{StrongRFProp}, \(A\times B\) is backwards-bounded and so by Theorem \ref{CStheorem}, is completely separable. If both are completely separable then for all \((a,b)\in A\times B\) there exists an \(n\in\mathbb{N}\) such that \(B_n(a)=B_n(b)=\varnothing\), and so \(B_n(a,b)=\varnothing\times f_2^{-n}(b)\cup f_1^{-n}(a)\times\varnothing=\varnothing\) and so \(A\times B\) is completely separable.

For the forward implication, we show the contrapositive. Assume without loss of generality that \(A\) is not completely separable and \(B\) is not backwards bounded. Then there exists \(a\in A\) and \(b\in B\) such that for all \(n\in\mathbb{N}\) \(B_n(a)\neq\varnothing\) and \(f_2^{-n}(b)\neq\varnothing\), and so \(B_n(a,b)\supseteq B_n(a)\times f_2^{-n}(b)\neq\varnothing\). Thus \(A\times B\) is not completely separable.
\end{proof}
\end{theorem}

\section{Concluding Remarks}

One application of the results in this paper is to investigate separation properties within varieties of monounary algebra. In \cite{jacobs1964lattice} a classification of all varieties of monounary algebras is given:

\begin{description}
    \item \(V_0\): The trivial variety, containing only the trivial monounary algebra \(T_1\), and the empty monounary algebra. This is defined by the equation \(x=y\) for all \(x,y\).
    \item \(V_k\) for \(k\in\mathbb{N}^+\): This contains all connected monounary algebras containing a cycle of length 1, such that for all \(x\) not in a cycle, \(f^{-k}(x)=\varnothing\). This is defined by the equation \(f^k(x)=f^k(y)\).
    \item \(V_{k,d}\) for \(k,d\in\mathbb{N}\), \(d\geq1\): This contains monounary algebras where each connected component contains a cycle whose length divides \(d\), and for all elements \(x\) not in a cycle, \(f^{-k}(x)=\varnothing\). This is defined by the equation \(f^{k+d}(x)=f^k(x)\).
    \item \(V_{0,0}\): The variety of all monounary algebras. This is defined be the equation \(x=x\).
\end{description}

We can note from these definitions that in any connected component of an algebra from \(V_{k,d}\) every vertex of the corresponding graph is either in a cycle of length at most \(d\), or is on a path of length at most \(k\) ending at such a cycle, and so \(B_{k+d+1}(x)=\varnothing\) for all \(x\), and thus the connected component is completely separable by Theorem \ref{CStheorem}. It can be similarly seen that any algebra in \(V_k\) is completely separable. Thus in every variety of monounary algebra other than the class of all monounary algebras, every algebra is completely separable (and thus both subalgebra separable and residually finite). It is then trivially true that in all but the full variety, we have \(A\times B\) is residually finite if and only if \(A\) and \(B\) are residually finite.

There are a few potential directions in which one could build on these results. Here we formulate a few questions of interest.

We saw in Example \ref{Subdirectproductexample} that we can have a residually finite subdirect product of two monounary algebras which are not residually finite. So perhaps it is possible to find some conditions on the components of the product and/or the construction of the subdirect product which ensure residual finiteness. Obviously one can also extend this question to subalgebra separability and complete separability.

\begin{question}
What are necessary and sufficient conditions for a subdirect product of monounary algebras to be residually finite?
\end{question}

Unary algebras are significantly different to monounary algebras. There are no obvious generalisations of the results from this paper that one could apply. However, the potential use of unary algebras to apply to more complex structures like semigroups, makes it a very intriguing topic for research.

\begin{question}
Can we find a criterion for residual finiteness in the more general class of unary algebras?
\end{question}

\textbf{Acknowledgments}

I would like to thank the EPSRC, AD Links Foundation, and the School of Mathematics and Statistics of the University of St Andrews for supporting me in my graduate research. I would also like to thank an anonymous referee for their comments on the paper in general, and specifically their insightful questions regarding varieties, which led to the application discussed at the beginning of the concluding remarks.

\bibliography{main}{}

\begin{thebibliography}{10}

\bibitem{ceccherini2010cellular}
T.~Ceccherini-Silberstein and M.~Coornaert.
\newblock {\em Cellular automata and groups}.
\newblock Springer Science \& Business Media, 2010.

\bibitem{evans1969some}
Trevor Evans.
\newblock Some connections between residual finiteness, finite embeddability
  and the word problem.
\newblock {\em Journal of the London Mathematical Society}, 2(1):399--403,
  1969.

\bibitem{ganyushkin2008classical}
O.~Ganyushkin and V.~Mazorchuk.
\newblock {\em Classical finite transformation semigroups: an introduction},
  volume~9.
\newblock Springer Science \& Business Media, 2008.

\bibitem{golubov1970finite}
E.~A. Golubov.
\newblock Finite divisibility in semigroups.
\newblock {\em Sib. Math. J.}, 11(6):920--931, 1970.

\bibitem{gray2009residual}
R.~Gray and N.~Ru{\v{s}}kuc.
\newblock On residual finiteness of direct products of algebraic systems.
\newblock {\em Monatsh. Math.}, 158(1):63--69, 2009.

\bibitem{jacobs1964lattice}
Eugene Jacobs and Robert Schwabauer.
\newblock The lattice of equational classes of algebras with one unary
  operation.
\newblock {\em The American Mathematical Monthly}, 71(2):151--155, 1964.

\bibitem{jakubikova2009monounary}
D~Jakub{\'\i}kov{\'a}-Studenovsk{\'a} and J~P{\'o}cs.
\newblock Monounary algebras.
\newblock {\em PJ {\v{S}}af{\'a}rik Univ., Ko{\v{s}}ice}, 2009.

\bibitem{katai2019fine}
Kamilla K{\'a}tai-Urb{\'a}n, Andr{\'a}s Pongr{\'a}cz, and Csaba Szab{\'o}.
\newblock The fine-and generative spectra of varieties of monounary algebras.
\newblock {\em Algebra universalis}, 80(2):22, 2019.

\bibitem{mal1983homomorphisms}
A.~I. Mal’cev.
\newblock On homomorphisms onto finite groups.
\newblock {\em Amer. Math. Soc. Transl. Ser. 2}, (119):67--79, 1983.

\bibitem{mayr2018finiteness}
P.~Mayr and N.~Ru{\v{s}}kuc.
\newblock Finiteness properties of direct products of algebraic structures.
\newblock {\em J. Algebra}, 494:167--187, 2018.

\bibitem{miller2020separability}
C.~Miller, G.~O'Reilly, M.~Quick, and N.~Ru{\v{s}}kuc.
\newblock On separability finiteness conditions in semigroups.
\newblock {\em arXiv preprint arXiv:2006.08499}, 2020.

\bibitem{MR0332990}
D.~J.~S. Robinson.
\newblock {\em Finiteness conditions and generalized soluble groups. {P}art 2}.
\newblock Springer-Verlag, New York-Berlin, 1972.

\bibitem{west2001introduction}
D.~West.
\newblock {\em Introduction to graph theory}, volume~2.
\newblock Prentice hall Upper Saddle River, 2001.

\end{thebibliography}

\bibliographystyle{plain}

\end{document}